\documentclass{article}
\usepackage{amsmath}

\begin{document}
\title{A Generalised Family of Ratio Product Estimator using Transformation Equation}
\author{Viplav Kumar Singh and *Rajesh Singh\\
Department of Statistics, Banaras Hindu University\\
           Varanasi-221005, India\\
           *Corresponding author}
\date{}
\maketitle
\begin{abstract}
In this paper we deal with the estimation of population variance of the study variable y using auxiliary information on variable x. A family of ratio and product-type estimators are proposed using suitable transformation on both random variable x(auxiliary variable) and y(study variable).Up to the first order of approximation the expression of mean square error and Bias term are obtained. An empirical study is carried out to illustrate the performance of the constructed estimator over others.

\end{abstract}

\section{Introduction}
~~~In survey sampling it is well recognized that the use of auxiliary information
 results in substantial gain in efficiency over the estimators which do not utilize such information.           In many manufacturing industries and pharmaceutical laboratories sometimes researchers are interested in the variation of their products. To measure the variations within the values of study variable y, the problem of estimating the population variance of study variable y received a considerable attention of the statistician in survey sampling including Isaki (1983), Singh and Singh (2001, 2003),  Jhajj et al. (2005), Kadliar and Cingi (2007), Singh et al. (2008),  Grover (2010), Singh et al. (2011) and Singh and Solanki (2012) have suggested improved estimator for estimation of ${S_{y}}^{2}$\\

~~~Let us consider a finite population $U=(U_{1},U_{2}....U_{N})$ having N units and let y and x are the study and auxiliary variable with means $\bar{y}$ and $\bar{x}$ respectively.The problem is to estimate the population variance $S_{Y}^{2}$  of the study variable y which is highly correlated with the auxiliary variable x .\\
Let a Simple Random Sample Without Replacement (SRSWOR) of size n is taken from a finite population of size N.Let $C_{y}=S_{y}/\bar{Y}$ and $C_{x}=S_{x}/\bar{X}$ be the population coefficient of variations of variables y and
x respectively,such that ${S_{y}}^{2}=\sum\limits_{i=1}^{N} (y_{i}-\bar{Y})^{2}/(N-1)$ and ${S_{x}}^{2}=\sum\limits_{i=1}^{N} (x_{i}-\bar{X})^{2}/(N-1)$ which is the Population variances for variables x and y.Let $\rho_{yx} $ be the population correlation coefficient between y and x,taking $c_{yx}=\rho_{yx} c_{y}c_{x}$.In this paper we use $\theta =1/n$ that is, we are neglecting the finite population correction factor in the coming computation of this article.Let $\mu_{pq}=\sum\limits_{i=1}^{N}(Y_{i}-\bar{Y})^{p}(X_{i}-\bar{X})^{q}$,Where ($\bar{Y}$,$\bar{X}$) denote the population mean of bivariate (y,x).Let $ \lambda_{pq}=\mu_{pq}/\mu_{20}^{p/2} \mu_{02}^{q/2}$ and taking $\beta_{2y}^{*}=(\beta_{2y}-1)$,$\beta_{2x}^{*}=(\beta_{2x}-1)$ and $\lambda_{22}^{*}=(\lambda_{22}-1)$.\\

\section{Some Known Estimators In Literature:}
~~~~~Some existing estimators in literature for estimating population variance  ${S_{y}}^{2}$ are -\\
1.Conventional unbiased estimator $\hat{S}_{y}^{2}=s_{y}^{2}$\\
\\
2.Usual ratio estimator $\hat{S}_{R}^{2}=s_{y}^{2}\left( \dfrac{S_{x}^{2}}{s_{x}^{2}}\right)$,defined by Isaki (1983).\\
\\
3.Usual regression estimator,$\hat{S}_{Reg}^{2}=s_{y}^{2}+b^{*}(S_{x}^{2}-s_{x}^{2})$,defined by Isaki (1983),
where b is the sample regression coefficient between $s_{x}^{2}$ and $s_{y}^{2}$.\\
\\
4.The Ratio-type estimator,$t_{k}=s_{y}^{2}\left[ \dfrac{(aS_{x}^{2}-b)}{\alpha(as_{x}^{2}-b)+(1-\alpha)(aS_{x}^{2}-b)}\right]$,given by Khoshnevisan et al. (2007), Singh et al.(2007).Where $\alpha$ is suitably chosen constant.\\
5.Another Ratio-type estimator, $t_{s}=s_{y}^{2}\left[2- \left( \dfrac{s_{x}^{2}}{S_{x}^{2}}\right)^{w} \right]$,given by Sahai and  Ray,where w is a suitably chosen constant.\\
\\
The expressions of biases of existing estimators ,up to the first order of approximation are given by-
 \begin{equation}
 Bias(\hat{S}_{y}^{2})=0
 \end{equation}
\begin{equation}
 Bias(\hat{S}_{R}^{2})=\theta {S}_{y}^{2}\left[ \beta_{2x}^{*}-\lambda_{22}^{*} \right] 
 \end{equation}
 \begin{equation}
 Bias(\hat{S}_{Reg}^{2})=0
 \end{equation}
 \begin{equation}
  Bias({t}_{k})={S}_{y}^{2}\left[ \alpha^{2}\gamma^{2}\beta_{2x}^{*}-\alpha \gamma \lambda_{22}^{*}\right] 
  \end{equation}
  Where,$\gamma=\dfrac{aS_{x}^{2}}{aS_{x}^{2}-b}$\\
  \begin{equation}
  Bias(t_{s})=S_{y}^{2}\theta \left[ \dfrac{w(w-1)}{2}\beta_{2x}^{*}-w \lambda_{22}^{*}\right] 
  \end{equation}
\\
To compare the efficiencies of these estimators, we require their mean square
errors. The expressions of the mean square errors (MSE), up to first order of
approximation, are as follows-
\begin{equation}
MSE(\hat{S}_{y}^{2})=\theta S_{y}^{4}\beta_{2y}^{*}
\end{equation}
\begin{equation}
MSE(\hat{S}_{R}^{2})=\theta S_{y}^{4}(\beta_{2y}^{*}+\beta_{2x}^{*}-2\lambda_{22}^{*})
\end{equation}
\begin{equation}
MSE(\hat{S}_{Reg}^{2})=\theta S_{y}^{4}\beta_{2y}^{*}\left[ 1-\dfrac{{\lambda_{22}^{*}}^{2}}{\beta_{2y}^{*}\beta_{2x}^{*}}\right] 
\end{equation}
\begin{equation}
MSE(t_{k})=\theta S_{y}^{4}\left[ \beta_{2y}^{*}+\alpha^{2} \gamma^{2} \beta_{2x}^{*}-2\alpha \gamma \lambda_{22}^{*}  \right] 
\end{equation}
where,  $\alpha(opt)=\dfrac{\lambda_{22}^{*}}{\gamma \beta_{2x}^{*} }$ such that $\gamma=\dfrac{aS_{x}^{2}}{aS_{x}^{2}-b}$ and a=1,b=1\\
\begin{equation}
MSE(t_{s})=\theta S_{y}^{4}\left[\beta_{2y}^{*}+w^{2} \beta_{2x}^{*}-2w \lambda_{22}^{*} \right] 
\end{equation}
where,$w(opt)=\dfrac{\lambda_{22}^{*}}{\beta_{2x}^{*}}$
\section{Proposed Class Of Estimator:}
~~~In this section we have proposed an estimator for estimating $S_{y}^{2}$ using transformation equation.Let us Suppose that a simple random sample of size n(less then N) is taken from the population without replacement and $s_{y}^{2}$ , $s_{x}^{2}$ are the sample variance of y and x respectively.\\
~~~~~~Let us suppose that,
\begin{align*}
 s_{ui}^{2}&=s_{yi}^{2}+a \\
~~~~ s_{vi}^{2}&=cs_{xi}^{2}+dS_{x}^{2} ~~~~~~~\forall i=1,2...N
\end{align*}
where a, c and d are known constants so that $s_{u}^{2}=s_{y}^{2}+a$ and $s_{v}^{2}=cs_{x}^{2}+dS_{x}^{2}$ are the sample variances of transformed variate u and v respectively and $S_{u}^{2}=S_{y}^{2}+a$ ,$S_{v}^{2}=(c+d)S_{x}^{2}$ are the corresponding population variances of u and v variate.Following Upadhyaya et al.(2011) and Singh et al.(2008) ,we have proposed a generalised class of estimator given as
\begin{equation}
t=\alpha_{1} s_{u}^{2}\left[ \dfrac{S_{v}^{2}}{\alpha s_{v}^{2}+(1-\alpha)S_{v}^{2}}\right]^{\beta}-a 
\end{equation}
~~~A large number of estimators may be identified as members of the proposed
family of estimators for suitable values of the scalars (a,c, d,$\alpha_{1},\alpha,\beta$) involved in suggested estimator-\\

~~~If $(a,c, d,\alpha_{1},\alpha,\beta)=(0,c,d,1,\alpha,0)$,then suggested estimator will become usual unbiased estimator. \\

~~~If $(a,c, d,\alpha_{1},\alpha,\beta)=(0,1,0,1,1,1)$,then suggested estimator will become usual ratio estimator.\\

~~~If $(a,c, d,\alpha_{1},\alpha,\beta)=(0,1,0,1,1,1)$ =$(0,1,0,1,1,-1)$,then suggested estimator will become product estimator.\\

To obtain the Bias and Mean square error of the proposed class of estimator we write-\\
\begin{equation*}
s_{y}^{2}=(1+e_{0})S_{y}^{2}
\end{equation*}
\begin{equation*}
s_{x}^{2}=(1+e_{0})S_{x}^{2}
\end{equation*}
Such that,\\
\noindent$E(e_{0})=E(e_{1})=0$,\\

\noindent$E(e_{0}^{2})=\theta \beta_{2y}^{*}$,\\

\noindent$E(e_{1}^{2})=\theta \beta_{2x}^{*}$,\\

\noindent$E(e_{0}e_{1})=\theta \lambda _{22}^{*}$\\ 
where,~~$\theta=\dfrac{1}{n}$,~~$\beta_{2y}=\lambda_{40}$ and $\beta _{2x}=\lambda_{04}$\\
\\
Expressing (11) in terms of $e_{i}'s$ we have,\\
\begin{equation}
t=\alpha_{1}(s_{y}^{2}+a)\left[ \dfrac{(c+d)S_{x}^{2}}{\alpha \left\lbrace c(1+e_{1})S_{x}^{2}+d S_{x}^{2}\right\rbrace +(1-\alpha)\left\lbrace(c+d)S_{x}^{2} \right\rbrace } \right]^{\beta}-a 
\end{equation}
After expanding right hand side of equation (12) we get -
\begin{equation}
t=\alpha_{1}(s_{y}^{2}+a)\left[ 1+A e_{1}\right] ^{-\beta}-a
\end{equation}
where~$ A=\dfrac{\alpha c}{c+d}$.\\

We assume that $A e_{1} \le 1$,so that equation (13)can be expandable in terms of power series up to the first order of approximation as,
\begin{equation}
t=\left[ \alpha_{1}S_{y}^{2}+\alpha_{1}e_{0}S_{y}^{2}+\alpha_{1} a\right] \left[ 1-\beta A e_{1}+\beta e_{1}^{2}\right]-a 
\end{equation}
Subtract $S_{y}^{2}$ from both sides of equation (14) and expanding it up to the first order of approximation we get,
\begin{equation}
t-S_{y}^{2}=\left[(\alpha_{1}-1)(S_{y}^{2}+a) +e_{0}\alpha_{1} S_{y}^{2}-e_{1}\alpha_{1}\beta A (S_{y}^{2}+a)
+e_{1}^{2} B \alpha_{1} (S_{y}^{2}+a)-\beta A S_{y}^{2} \alpha_{1}e_{0}e_{1}\right] 
\end{equation}
Taking expectations both sides,we get bias of 't' up to the first order of approximation as,\\
\begin{equation}
Bias(t)=\left[ (\alpha_{1}-1)(S_{y}^{2}+a)+B\alpha_{1}(S_{y}^{2}+a)\theta \beta_{2x}^{*}-\beta A S_{y}^{2} \alpha_{1} \theta \lambda_{22}^{*}\right] 
\end{equation}
Squaring equation (15) and retaining terms up to the first order of approximation, we get the MSE of the  estimator t as,
\begin{equation}
MSE(t)=\left[(\alpha_{1}-1)^{2}(S_{y}^{2}+a)^{2}+\alpha_{1}^{2}Q_{1}-2\alpha_{1}Q_{2}\right]
\end{equation}
where,
\begin{align*}
Q_{1}&=\theta\left[ S_{y}^{4}\beta_{2y}^{*}+\beta_{2x}^{*}\left\lbrace\beta^{2}A^{2}(S_{y}^{2}+a)^{2}+2B(S_{y}^{2}+a)^{2} \right\rbrace-4\beta A S_{y}^{2}(S_{y}^{2}+a) \lambda_{22}^{*} \right] \\
\\
Q_{2}&=\theta \left[ B(S_{y}^{2}+a)^{2}\beta_{2x}^{*}-S_{y}^{2}(S_{y}^{2}+a) \beta A \lambda_{22}^{*} \right]\\
\end{align*}
and $\alpha_{1}(opt)=\left[\dfrac{Q_{2}+(S_{y}^{2}+a)^{2}}{Q_{1}+(S_{y}^{2}+a)^{2}} \right] .$\\
After putting $\beta=1,\alpha=1$,$a=C_{x},b=C_{x},d=0.9742$ in the expression of $Q_{1}and Q_{2}$ and then putting these value of $Q_{1},Q_{2}$and $\alpha_{1}$ in equation (17),our proposed estimator t tends to regression estimator.
\section{Numerical Illustration}
To have a rough idea about the gain in efficiency of the proposed and existing estimators, defined under
the situations when prior information of population variance of auxiliary variable
is available, we take the same empirical data as considered by Shabbir and Gupta
(2007) and Kadilar and Cingi (2007). \\
Consider 104 villages of the East Anatolia
Region in Turkey.The variables are-\\
\textbf{Data statistics:}\\
y : level of apple production (1 unit = 100 tones)\\
x : number of apple trees (1 unit = 100 trees).\\
\\
~~~~The required value of population parameters are:
\begin{center}
$N=104,S_{y}=11.6694,S_{x}= 23029.072,C_{y}=1.866,C_{x}=1.653,$\\
$\rho_{yx}=0.865, C_{yx}=2.668, \beta_{2y}=16.523, \beta_{2x}=17.516, \lambda_{22}=14.398$
\end{center}
\noindent Percent Relative Efficiency (PRE) of an estimator is given by-\\
\begin{equation*}
PRE(.)=\dfrac{var(\hat{S}_{y^{2}})}{MSE(.)}
\end{equation*}\\
\newpage
The Mean Square Error and percent relative efficiencies of various estimators are given in the Table 1
\begin{center}
\begin{tabular}{p{2 cm}p{3 cm}p{2 cm}p{2 cm}}
\hline
\textbf{Estimators} &\textbf{Mean Square Error}  &\textbf{Percent Relative Efficiencies} \\ 
\hline
 $\hat{S}_{y}^{2}$& 14395.4 & 100.00 \\ 
\\ 
 $\hat{S}_{R}^{2}$& 4862.145 & 296.071 \\ 
\\
 $\hat{S}_{Reg}^{2}$&4316.267  &333.515  \\ 
 \\
 $t_k$& 4316.267  & 333.515 \\ 
\\
 $t_s$& 4316.267 & 333.515 \\ 
\\
$t$ & 4316.258 & 333.515 \\ 
\hline
\end{tabular}
\end{center}

\section{Conclusion }
From table 1 we observed that the proposed estimator t performs better than the estimator $\hat{S}_{y}^{2}$ and Isaki(1983) estimator $\hat{S}_{R}^{2}$ .we also observed that proposed estimator under optimum conditions performs equally efficient as regression estimator $\hat{S}_{Reg}^{2}$. 
\bibliographystyle{plain}

\begin{thebibliography}{1}

\bibitem{}
Grover, L. K . (2010): A correction note on improvement in variance estimation using auxiliary information. Commun. Stat. Theo. Meth. 39:753–764.
\bibitem{}
Isaki, C.T.(1983): Variance estimation using auxiliary information.  Journal of American Statistical Association 78, 117–123.
\bibitem{}
Jhajj, H .S., Sharma, M. K. and Grover, L. K. (2005) : An efficient class of chain estimators of population variance under sub-sampling scheme. J. Japan Stat. Soc., 35(2), 273-286.
\bibitem{}
Kadilar, C. and Cingi, H. (2006) :  Improvement in variance estimation using auxiliary information. Hacettepe Journal of Mathematics and Statistics 35 (1), 111–115.

\bibitem{}
Singh, R.,Cauhan, P., Sawan, N. and Smarandache,F. (2007):  Auxiliary information and a priori values in construction of improved estimators. Renaissance High press.
\bibitem{}
Singh, H.P. and  Singh,  R. (2001):  Improved ratio-type estimator for variance using auxiliary information. J Indian Soc Agric Stat 54(3):276–287.

\bibitem{}
Singh,H.P. and Singh, R. (2003):  Estimation of variance through regression approach in two phase sampling. Aligarh Journal of Statistics, 23, 13-30.
\bibitem{}
Singh, H. P. and Solanki, R.S. (2012):  A new procedure for variance estimation in simple random sampling using auxiliary information. Stat. Paps.  DOI 10.1007/s00362-012-0445-2.
\bibitem{}
Singh, R., Chauhan, P., Sawan, N. and Smarandache,F.(2008):   Almost Unbiased Ratio and Product Type Estimator of Finite Population Variance Using the Knowledge of Kurtosis of an Auxiliary Variable in Sample Surveys. Octogon Mathematical Journal, Vol. 16, No. 1, 123-130.
\bibitem{}
Singh, R., Chauhan, P., Sawan, N. and Smarandache, F.(2011): Improved Exponential Estimator for Population Variance Using Two Auxiliary Variables. Italian Jour. Of Pure and Appld. Math., 28, 103-110.
\bibitem{}
Upadhyaya L.N Singh H.P Chatterjee S. Yadav R.(2011):A generalized family of transformed ratio-product
estimators in sample surveys ,Model Assisted Statistics and ApplicationsIOS Press 	Volume 6, Number 2 / 2011,Pages 137-150.

\end{thebibliography}

\end{document}